\documentclass{article}
\usepackage{amsmath}

\newcommand{\remove}[1]{}
\usepackage{tikz}
\usepackage{hyperref}
\usepackage{forest}

\begin{document}
\title{A High Quartets Distance Construction}
\author{Benny Chor\footnote{School of Computer Science, Tel Aviv University. {\tt benny@cs.tau.ac.il}}
\and P\'eter L. Erd\H{o}s\footnote{MTA  A. R\'enyi Institute of Mathematics, Budapest. {\tt erdos.peter@renyii.mta.hu}}
\and Yonatan Komornik\footnote{School of Computer Science, Tel Aviv University. {\tt yoniko@gmail.com}} }
\maketitle
\begin{abstract}
Given two binary trees on $N$ labeled leaves, the {\em quartet distance} between the trees
 is the number of disagreeing quartets.
By permuting the leaves at random, the expected quartets distance between the two trees is $\frac{2}{3}\binom{N}{4}$.
However, no strongly explicit construction reaching this bound asymptotically was known.

We consider complete, balanced binary trees on $N=2^n$ leaves,
labeled by $n$ long bit sequences. Ordering the leaves in one tree by the prefix order, and in the other
tree by the suffix order, we show that the resulting quartet distance is $\left(\frac{2}{3} + o(1)\right)\binom{N}{4}$,
and it always exceeds the $\frac{2}{3}\binom{N}{4}$ bound.
\end{abstract}

\section{Background}
Given a set of taxa (a group of related biological species), the goal of phylogeny reconstruction is to build a tree which best represents the course of evolution for this set over time. The leaves of the tree are labeled with the given, extant taxa. Internal nodes correspond to hypothesized, extinct taxa. There are numerous phylogeny reconstruction approaches \cite{SempleSteel}. One approach of interest is building unrooted, resolved (or binary) trees from quartets, where a quartet is an unrooted tree on 4 leaves. We note that for a given set of 4 leaves there are 3 quartets topologies. The input is a set of (possibly weighted) quartets, and the goal is to build a tree which would agree with the maximum number of input quartets (maximum weighted sum, correspondingly)  \cite{BC, BJ, SR, SV}. It is known that this problem is computationally hard \cite{Steel}.

Various combinatorial problems related to quartets have also been studied extensively. In this paper, we are especially interested in the {\em quartet distance} problem \cite{EMM}. Let $T_1, T_2$ be two resolved (binary) trees on the same set of $N$ labeled leaves. Every set of the same 4 leaves induces two quartets, one in $T_1$ and the other in $T_2$. The topologies of the two quartets could either agree or disagree. The quartet distance between $T_1, T_2$ is the number of disagreeing  quartets. Notice that the identity of a quartet in a given
binary tree is well defined, regardless of the placement of the root. Thus the quartet distance between $T_1, T_2$ is invariant under different
rootings of $T_1, T_2$, and under making one or both trees unrooted. We remark that
there are efficient algorithms to compute the quartet distance of two trees. The most efficient one, by Brodal, Fagerberg, and Pedersen, runs in $O(N\log N)$ time \cite{Brodal}.

Bandelt and Dress \cite{BD} conjectured that the maximum quartet distance between any two resolved (binary) trees on $N$ leaves is at most $\left(\frac{2}{3} + o(1)\right)\cdot\binom{N}{4}$. Taking  two  binary trees $T_1, T_2$ on the same set of $N$ leaves,  and assigning labels to the leaves at random, the probability that any quartet will agree equals exactly $1/3$. This implies that the expected value of the quartet distance is exactly  $\frac{2}{3}\cdot \binom{N}{4}$.  This simple probabilistic argument can be de-randomized using standard de-randomization methods. We will further refer to the result of such de-randomization in the context of our work in Section~\ref{conc}.

Alon, Snir, and Yuster \cite{ASY} showed that the random labeling method implies the existence of trees with quartet distance strictly greater than $\frac{2}{3}\cdot \binom{N}{4}$. They also proved a $\frac{9}{10}\cdot \binom{N}{4}$ asymptotic upper bound on the quartet distance. Finally, using the technique of flag algebra, Alon, Naves, and Sudakov \cite{ANS} have shown a $\left(0.69 +  o(1)\right) \cdot\binom{N}{4}$ upper bound on the normalized quartet distance (for large enough number of leaves, $N$).

No strongly explicit construction attaining the $\frac{2}{3}\cdot \binom{N}{4}$ lower bound asymptotically is known (the notions of 
explicit and strongly explicit constructions are defined and discussed in Section~\ref{conc}). We consider complete, balanced binary trees on $N=2^n$ leaves, labeled by $n$ long bit sequences. Ordering the leaves in one tree by the {\bf prefix} (or lexicographic) order, and in the other tree by the {\bf suffix} (or co-lexicographic) order, we show that the resulting quartet distance is $\left(\frac{2}{3} + o(1)\right)\cdot \binom{N}{4}$, and furthermore, the distance exceeds the $\frac{2}{3}\cdot\binom{N}{4}$ bound for all $N$. An important part of our proof  is counting 
the number of binary strings whose longest common prefixes (or suffixes) are of given lengths.

\section{High Level View}

Denote by $Pref_n$ the complete, balanced binary tree with leaves labeled by $\{0,1\}^n$ and ordered by prefix (or lexicographic) order, and by $Suff_n$ the complete, balanced binary tree on the same set of leaves, ordered by suffix (or co-lexicographic) order. Consider an ordered 4-tuple of distinct binary sequences $(x_0,x_1,x_2,x_3)$, $x_i\in \{0,1\}^n$ (these are the labels of leaves in our two trees). For every pair of indices $0\leq i< j\leq 3$, let $P_{i,j}(x_0,x_1,x_2,x_3)$ be the event ``the common prefix of $x_i, x_j$ is not shorter than the other five common prefixes''. Likewise, we define the event $S_{i,j}(x_0,x_1,x_2,x_3)$, referring to suffixes. For sake of brevity, we will use $P_{i,j}, S_{i,j}$  to denote $P_{i,j}(x_0,x_1,x_2,x_3), S_{i,j}(x_0,x_1,x_2,x_3)$, correspondingly.

\begin{figure}[hbt]
\label{complete4}
\scalebox{.9}{
\begin{forest}
[$\bullet$
    [$\bullet$
        [$\bullet$, s sep=1mm,
            [{000}]
            [{001}]
        ]
        [$\bullet$
            [{010}]
            [{011}]
        ]
    ]
    [$\bullet$
        [$\bullet$
            [{100}]
            [{101}]
        ]
        [$\bullet$
            [{110}]
            [{111}]
        ]
    ]
]
\end{forest}
}
\qquad
\scalebox{.9}{
\begin{forest}
[$\bullet$
    [$\bullet$
        [$\bullet$
            [{000}]
            [{100}]
        ]
                [$\bullet$
            [{010}]
            [{110}]
        ]
    ]
    [$\bullet$
        [$\bullet$
            [{001}]
            [{101}]
        ]
        [$\bullet$
            [{011}]
            [{111}]
        ]
    ]
]
\end{forest}
}
\caption{The complete, balanced binary tree for strings of length $n=3$, with labels in prefix (lexicographic) order on the left, and
suffix (co-lexicographic) order on the right.}
\end{figure}
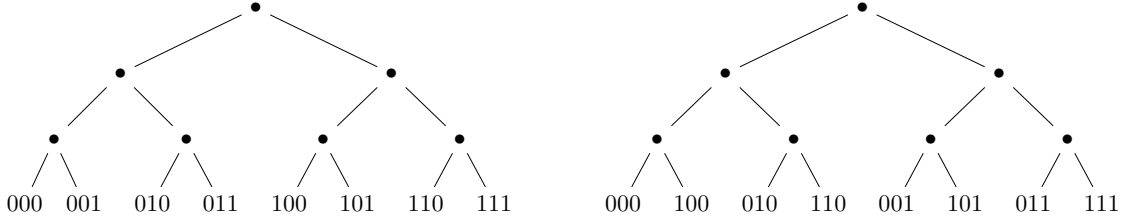

There are some obvious relations among the  $P_{i,j}$ or the $S_{i,j}$. For example $P_{0,1}$,  $P_{0,2}$,  $P_{0,3}$ are mutually exclusive. More generally, any pair $P_{i_1,j_1}$, $P_{i_2,j_2}$ sharing exactly one subscript ($i_1=i_2$ or $j_1=j_2$) is mutually exclusive. Note, however, that {\em e.g.} $P_{0,1}$ and $P_{2,3}$ are {\em not} mutually exclusive. Clearly, the number of ordered binary sequences satisfying $P_{i,j}, S_{i,j}$ is the same for all choices of indices $i<j$.

To determine the quartet distance between our two trees, we will exactly compute the number of length $n$ sequences satisfying various combinations of these events, such as $P_{0,1} \cap P_{2,3}$, $P_{0,1} \cap S_{0,1}$, $P_{0,1} \cap P_{2,3} \cap S_{0,1}$, and $P_{0,1} \cap P_{2,3} \cap S_{0,1} \cap S_{2,3}$. These, in turn, will enable the derivation of the exact and asymptotic quartet distance between the ``suffix order'' and the ``prefix order'' binary sequences' trees, using a simple inclusion-exclusion argument.

\section{Sequence Counts of Specific Events}
For each event we will count the number of four tuples of different {\em ordered} sequences,
$(x_0, x_1, x_2, x_3)$
satisfying it, using simple properties of prefixes and suffixes of $n$ bit long binary sequences. The lengths of common prefixes
and suffixes of any pair of binary sequences remains invariant by xoring the sequences to any one sequence
(namely computing the bit-wise XOR of the sequences). By xoring the four sequences
to $x_0$, we can thus assume without loss of generality that $x_0$ is the all $0$ sequence, while $x_1, x_2, x_3$ are three
uniformly distributed sequences that are non zero and distinct.

\subsection{$P_{0,1} \cap S_{2,3}$}
Let us denote the length of the longest common prefix of $x_0,x_1$ by $\ell, (\ell\leq n-1)$, and the length of the longest common suffix of  $x_2,x_3$ by $k, (k\leq n-1)$. For $P_{0,1}$, $\ell\geq 1$ should hold, and for $S_{2,3}$, $k\geq 1$ should hold. We treat separately the following three cases:
\begin{enumerate}
\item[(1)] $\ell+k+2\leq n$  (the $\ell$ long prefix plus one bit buffer zone, and the $k$ long suffix plus one bit buffer
zone, do not overlap). Note that since $1\leq k$, the value $\ell$ is bounded by $\ell\leq n-3$.
\item[(2)] $\ell+k+1 = n$
\item[(3)] $\ell + k \geq n$
\end{enumerate}

We start by analyzing case (1). There is no overlap between the $\ell+1$ long prefixes  and the $k+1$ long  suffixes. This will enable us to analyze the number of possible prefixes and possible suffixes for $x_1, x_2, x_3$ separately, thereby facilitatating the counting. Let us start with the prefixes.  Given that $x_0=0^n$,  as the longest common prefix of $x_0,x_1$ is of length  $\ell$, the $\ell+1$ long prefix of $x_1$ must be $0^{\ell}1$. The $\ell+1$ long prefix of $x_2$  must differ from both $0^{\ell+1}$ and $0^{\ell}1$ for the event $P_{0,1}$ to hold. Thus, there are $2^{\ell+1}-2=2(2^\ell-1)$ possibilities for choosing the $\ell+1$ long prefix of $x_2$. By a similar argument,
there are $2^{\ell+1}-3$ possibilities for choosing the $\ell+1$ long prefix of $x_3$. So, given that $x_0=0^n$, the number of possibilities for choosing the $\ell+1$ long prefixes of $x_1, x_2, x_3$ is
$2\cdot (2^\ell-1)\cdot (2^{\ell+1}-3)=2^{2\ell+2}-5\cdot 2^{\ell+1}+6$.

Let us now turn to the suffixes. Let $b_0b_1\ldots b_{k-1}b_k\in\{0,1\}^k$ denote the $k+1$ long suffix of $x_2$. This determines uniquely the  $k+1$ long suffix of $x_3$, which equals   $\overline{b}_0 b_1\ldots b_{k-1}b_k\in\{0,1\}^k$.
For $S_{2,3}$ to hold, both should differ from the $k+1$ long suffix of $x_0$, which equals $0^{k+1}$. In particular, the $k+1$ long suffix of $x_2$ must differ from both $0^{k+1}$ and $10^k$. This leaves $2^{k+1}-2=2(2^k-1)$ possibilities for choosing the $k+1$ long suffix of $x_2$, and then $2^{k+1}-3$ possibilities for choosing the $k+1$ long suffix of $x_1$. So, given that $x_0=0^n$, the number of possibilities for choosing the $k+1$ long suffixes of $x_1, x_2, x_3$ is  $2\cdot (2^k-1)\cdot (2^{k+1}-3)=2^{2k+2}-5\cdot 2^{k+1}+6$.

Finally, each of $x_1, x_2, x_3$ has $n-\ell-k-2$ ``free bits'' in the middle, not overlapping neither the prefix nor the suffix. These can vary over all possibilities, independently of each other. The total number of possibilities for the free bits of the three sequences is thus $2^{3n-3\ell-3k-6}$, and the total number of possibilities for all of $x_1, x_2, x_3$, given that $x_0=0^n$, is
\begin{equation*}
\begin{aligned}
&2^{3n-3\ell-3k-6}\cdot \left( 2^{2\ell+2}-5\cdot 2^{\ell+1}+6 \right)\cdot\left(2^{2k+2}-5\cdot 2^{k+1}+6\right)\\
=&2^{3n-4}\cdot \left( 2^{-\ell+1}-5\cdot 2^{-2\ell}+3\cdot2^{-3\ell} \right)\cdot\left(2^{-k+1}-5\cdot 2^{-2k}+3\cdot2^{-3k}\right)\\
\end{aligned}
\end{equation*}
Summing over $\ell$ and $k$ in the relevant range, we get
$$
2^{3n-4}\cdot\sum_{\ell=1}^{n-3}\left( 2^{-\ell+1}-5\cdot 2^{-2\ell}+3\cdot2^{-3\ell} \right)\cdot\sum_{k=1}^{n-\ell-2}
\left(2^{-k+1}-5\cdot 2^{-2k}+3\cdot2^{-3k}\right)\ (*).$$\\
Employing a symbolic algebra package (specifically, Maple) to this sum, we get
$${\frac {16}{441}}\cdot 2^{3n}-n\cdot2^{2n} + 5\cdot 2^{2n}-{\frac {25}{3}}\cdot n\cdot{2}^{n}+{\frac {95}{9}}\cdot {2}^{n}
-{\frac {36}{7}}\cdot n-{\frac {764}{49}}\ .$$
This is the number of ordered quartets with $x_0=0^n$, satisfying case (1) of
 $P_{0,1} \cap S_{2,3}$.\\

\bigskip
Let us turn to case (2), where $\ell+k+1=n$, which means that the $\ell$ bits long prefix and the $k$ bits long suffix do not overlap, and have one ``buffer bit'', which separates them. For the
event $S_{2,3}$ to occur and the longest common suffix of $x_2$, $x_3$ to be of length
$k$, the $k$ bits long suffix of $x_2$ must differ from the $k$ bits long suffixes of both
$x_0$ and $x_1$.

Given $P_{0,1}, x_0 = 0^n$, and $\ell$ being the length of the
longest common prefix of $x_0, x_1$, the $\ell+1$ bits long prefix of $x_1$ is $0^\ell 1$. Since $\ell+k+1=n$, the last bit of the $\ell+1$ bits long prefix of $x_1$ is also the first bit of its  $k+1$ bits long suffix. So this suffix differs from the $k+1$ bits long suffix of $x_0$ (which equals $0^{k+1}$).

There are $2^{k}$ possible settings of the $k$ rightmost bits of $x_1$. We treat separately
the case (a) where these bits are $0^k$, and the case (b) where they differ from $0^k$. In case
(a), neither $10^k$ nor $0^{k+1}$ can serve as the $k+1$ bits suffix of $x_2$, but any other sequence can. There are $2^{k+1}-2$ such possibilities. Given the $k+1$ bits suffix of $x_2$, the $k+1$ bits suffix of $x_3$ is completely determined (it differs from $x_2$ in the
buffer zone bit, and agrees with it in the other $k$ bits). The $\ell$ bits long prefix of $x_2$ and of
$x_3$ could be any two sequences, other than $0^{\ell}$. So in case (a) the overall number of possibilities for $x_1, x_2, x_3$ is
$$1\cdot (2^{k+1}-2)\cdot (2^{\ell}-1)^2=2\cdot(2^k-1)\cdot(2^{\ell}-1)^2\ .$$

In case (b), there are $2^k-1$ possibilities for the $k$ bits long suffix of $x_1$. The $k$ rightmost bits of $x_2$ must differ from both the $k$ rightmost bits of $x_1$, and from
$0^k$. Thus, there are $2^k-2$ possibilities for the $k$ bits long suffix of $x_2$, and $2$
possibilities for the buffer zone bit of $x_2$. Overall,  this leaves $2\cdot(2^{k}-2)$ possibilities for the $k+1$ bits suffix of $x_2$, which completely determine the $k+1$ bits suffix of $x_3$.
Like case (a), the $\ell$ bits long prefix of $x_2$ and of $x_3$ have $2^{\ell}-1$ possibilities each.
So in case (b) the overall number of possibilities for $x_1, x_2, x_3$, for a given value
of $k$ and $\ell$, is
$$2\cdot(2^k-1)\cdot(2^k-2)\cdot(2^{\ell}-1)^2\ .$$

Summing the numbers in cases (a) and (b), we get
\begin{equation*}
\begin{aligned}
&2\cdot(2^k-1)\cdot(2^{\ell}-1)^2+2\cdot(2^k-1)\cdot(2^k-2)\cdot(2^{\ell}-1)^2\\
=&\ (2^k-1)\cdot(2^{\ell}-1)^2\cdot\left(2+2\cdot(2^k-2)\right)\\
=&\ 2\cdot(2^k-1)^2\cdot(2^{\ell}-1)^2\ .
\end{aligned}
\end{equation*}
In case (2) $\ell+k+1=n$, so $\ell=n-k-1$. Furthermore, $k,\ell\geq 1$, so $k$
is in the range $1\leq k\leq n-2$. Summing over all values of $k$, we get that
the number of ordered quartets with $x_0=0^n$, satisfying case (2) of
 $P_{0,1} \cap S_{2,3}$, equals\\
\begin{equation*}
\begin{aligned}
&\sum_{k=1}^{n-2}2\cdot\left(2^k-1\right)^2\cdot\left(2^{n-k-1}-1\right)^2\\
&=\frac{1}{2}\cdot n\cdot 2^{2n}-\frac{8}{3}\cdot 2^{2n}+4\cdot n\cdot 2^{n}-4\cdot 2^n+2\cdot n+\frac{20}{3}\ .
\end{aligned}
\end{equation*}

\bigskip
We will now turn to case (3), where $n\leq \ell+k$, so there is no buffer bit between the $\ell$ bits
long prefix and the $k$ bits long suffix, and if $n < \ell+k$, they even overlap.
Again, we assume that $x_0 = 0^n$, thus the $\ell + 1$ leftmost bits of $x_1$ are $0^\ell 1$.
We then have $2^{n-\ell-1}$ ways to choose $x_1$'s suffix. Since $n-\ell-1 < k$, it is guaranteed that even if $x_0$ and $x_1$ shared $n-\ell-1$ suffix bits, their common suffix won't be longer than $k$.

Now, given $x_0$ and $x_1$, we want to determine the number of
possibilities for $x_2$ and $x_3$. Note that the $n-k-1$ bit of $x_2$ and $x_3$ must differ (otherwise
the length of the common suffix would be greater than $k$.

Consider the $(k+\ell)-n$ bits of $x_2, x_3$, where the $\ell$ long prefix
and $k$ long suffix overlap (if $k+\ell = n$, this overlap is empty). These bits are part of the $k$ long suffix,
shared by $x_2$ and $x_3$. Let us consider the 2 following sub cases:
\begin{enumerate}
\item[(i)] The $(k+\ell)-n$ bits of $x_2, x_3$ equal $0^{k+\ell-n}$.

In this case, the $n-\ell$ rightmost bits of of $x_2, x_3$ must differ from the $n-\ell$ rightmost bits of
$x_0$ (which are all $0$) and of $x_1$ (which are not all $0$). The number of possibilities is thus
$2^{n-\ell}-2$. Suppose, without loss of generality, that the $n-k-1$ bit of $x_2$ equals $0$. The
$n-k-1$ long prefix of $x_2$ must differ from $0^{n-k-1}$ (otherwise it would share an $\ell$ long prefix with
both $x_0$ and $x_1$). There are $2^{n-k-1}-1$ possibilities for this prefix. There is no such restrictions on the
$n-k-1$ bit long prefix of $x_3$, so there are $2^{n-k-1}$ possibilities for it. Overall,
the number of possible sequences in case (i) is $2\cdot 2^{n-\ell-1}\cdot (2^{n-\ell}-2)\cdot (2^{n-k-1}-1)\cdot 2^{n-k-1}$,
where the leading $2$ accounts for the cases where either the $n-k-1$ bit of $x_2$  or that bit of $x_3$
equals $0$.

\item[(ii)]  The $(k+\ell)-n$ bits of $x_2, x_3$ differ from $0^{k+\ell-n}$.

There are $2^{(k+\ell)-n}-1$ ways to determine these $(k+\ell)-n$ bits of $x_2, x_3$.
And there are $2^{n-\ell}$ ways to determine the $n-\ell$ bit long suffix of $x_2, x_3$.
Suppose, without loss of generality, that the $n-k-1$ bit of $x_2$ equals $0$. There are no additional restrictions on the $n-k-1$ long prefix of $x_2$,
so there are $2^{n-k-1}$ possibilities for this prefix. There are exactly that many possibilities for the $n-k-1$
long prefix of $x_3$.
Overall,
the number of possible sequences in case (ii) is $2\cdot (2^{(k+\ell)-n}-1)\cdot 2^{n-\ell-1}\cdot 2^{n-\ell} \cdot 2^{n-k-1}\cdot 2^{n-k-1}$,
where the leading $2$ accounts for the cases where either the $n-k-1$ bit of $x_2$  or that bit of $x_3$
equals $0$.

\end{enumerate}

Summing up cases (i, ii), we get that the number of possibilities for $x_1, x_2, x_3$ equals
\begin{equation*}
\begin{aligned}
& 2\cdot 2^{n-\ell-1}\cdot (2^{n-\ell}-2)\cdot (2^{n-k-1}-1)\cdot 2^{n-k-1}\\
+& 2\cdot (2^{(k+\ell)-n}-1)\cdot 2^{n-\ell-1}\cdot 2^{n-\ell} \cdot 2^{n-k-1}\cdot 2^{n-k-1}\\
\end{aligned}
\end{equation*}
Summing over values of $k$ and $\ell$, satisfying $n\leq k+\ell$, we get

\begin{equation*}
\begin{aligned}
\sum_{\ell=1}^{n-1}\sum_{k=n-\ell}^{n-1}&\left(2\cdot 2^{n-\ell-1}\cdot (2^{n-\ell}-2)\cdot (2^{n-k-1}-1)\cdot 2^{n-k-1}\right.
\\
&+ \left. 2\cdot (2^{(k+\ell)-n}-1)\cdot 2^{n-\ell-1}\cdot 2^{n-\ell} \cdot 2^{n-k-1}\cdot 2^{n-k-1}\right)\\
=& \frac{1}{2}\cdot n\cdot 2^{2n}-{\frac {7}{3}}\cdot 2^{2n}+2\cdot n \cdot 2^n +2^n +\frac{4}{3}\ .
\end{aligned}
\end{equation*}
Summing the contributions from cases (1), (2), and (3), we conclude that the number of ordered quartets
with $x_0=0^n$ in $P_{0,1}\cap S_{2,3}$ equals
$$
 \frac{16}{441}\cdot 2^{3n}-  \frac{7}{3}\cdot n\cdot2^{n}+\frac{68}{9}\cdot 2^{n}-\frac{22}{7}\cdot n - \frac{372}{49}\ .
$$
Note that the $\theta(n\cdot2^{2n}),\ \theta(2^{2n})$ terms were cancelled.

\subsection{$P_{0,1} \cap S_{0,1}$}

We denote the length of the longest common prefix of $x_0,x_1$ by $\ell\  (\ell\leq n-1)$, and the length of the longest common suffix of  $x_0,x_1$ by $k\  (k\leq n-\ell-1)$. For $P_{0,1}$, $\ell\geq 1$ should hold, and for $S_{0,1}$, $k\geq 1$ should hold.
Note that in this case, the locations of the longest common suffix and the longest common prefix cannot intersect .
We treat separately the following two cases:
\begin{enumerate}
\item[(1)] $\ell+k+2\leq n$  (the $\ell$ long prefix plus one bit buffer zone, and the $k$ long suffix plus one bit buffer
zone, do not overlap). Since $1\leq k$, $\ell$ is bounded by $\ell\leq n-3$.
\item[(2)] $\ell+k+1 = n$.
\end{enumerate}
Note that $\ell + k < n$ must hold, for otherwise we would have $x_0=x_1$.
Given that $x_0$ is $0^n$, it is then clear that $x_1$'s $\ell+1$ long prefix is $0^\ell 1$ and its $k+1$ long suffix is $1 0^k$.\\

In case (1), $\ell + k + 1 < n$, and $x_1$ has the form $x_1 = 0^\ell 1 x 1 0^k$  where $x \in \{0,1\}^{n-k-\ell-2}$.
$x$ can be chosen with no constrains from $\{0,1\}^{n-k-\ell-2}$, so there are $2^{n-k-\ell-2}$ ways to choose $x$.
There are $2^{\ell+1}-2$ ways to choose the $\ell+1$ long prefix of $x_2$ (it must differ from the $\ell+1$ long prefix of
$x_0$ and $x_1$), and $2^{\ell+1}-3$ ways to choose the $\ell+1$ long prefix of $x_3$ (it must differ from the $\ell+1$ long prefixes of
$x_0$, $x_1$, and $x_2$). In a similar manner, there are $2^{k+1}-2$ ways to choose the $k+1$ long suffix of $x_2$, and
$2^{k+1}-3$ ways to choose the $k+1$ long suffix of $x_3$. Finally, the remaining
$n-k-\ell-2$ bits of the buffer zone in each of $x_2, x_3$ can be chosen freely. All by all, the number of possibilities of case (1)
for given values of $\ell$ and $k$ is $$2^{3(n-k-\ell-2)}\cdot (2^{\ell+1}-2)\cdot (2^{\ell+1}-3) \cdot (2^{k+1}-2)\cdot (2^{k+1}-3)\ .$$
We remark that this expression is the same as the one derived for case (1) of $P_{0,1} \cap S_{2,3}$.

In case (2), $\ell + k + 1 = n$, and $x_1$ has the form $x_1 = 0^\ell10^k$, so it is completely determined. Unlike
case (1), the $\ell+1$ long prefix and $k+1$ long suffix overlap, which makes the treatment slightly more involved.
We therefore partition case (2)  into two subcases: (i) $x_2$ and $x_3$'s common suffix length is shorter than $k$, and (ii) $x_2$ and $x_3$'s common suffix length is exactly $k$.
For case (i) we can still choose $x_2$ and $x_3$'s prefixes as we have done in (1), namely there are $(2^{\ell+1}-2) \cdot (2^{\ell+1}-3)$ ways to choose them. Given the $\ell+1$ long suffixes, both $x_2$ and $x_3$ still got $n - \ell - 1 = k$ bits that are not yet determined. Since the length of their shared suffix is shorter than $k$, the two choices must be different from each other, and from $0^k$. So there are $(2^{k}-1) \cdot (2^{k}-2)$ ways to choose the remaining $k$ bits. All by all, the number of possibilities in subcase (i) is
$(2^{\ell+1}-2) \cdot (2^{\ell+1}-3)\cdot (2^{k}-1) \cdot (2^{k}-2)$ possibilities.

For subcase (ii), $x_2$ and $x_3$'s $k$ long suffixes are the same, but the $(k+1)$th bits (from the right) are different.
The $k+1$ long suffixes must be different from both $0^{k+1}$ and $10^k$.
This leaves $2^{k+1}-2$ choices for $x_2$'s $k+1$ long suffix, and determines $x_3$'s $k+1$ long suffix. Now the $\ell$ long prefixes of both can be chosen freely, as long as they both are not $0^\ell$. So there are $(2^\ell-1)^2$ ways to choose the prefixes for $x_2$ and $x_3$. In subcase (2)(ii) there are $(2^\ell-1)^2\cdot (2^{k+1}-2)$ possibilities.
All by all, the number of possibilities in case (2) for given values of $\ell$ and $k$ is
$$(2^{\ell+1}-2) \cdot (2^{\ell+1}-3)\cdot (2^{k}-1) \cdot (2^{k}-2)+(2^\ell-1)^2\cdot (2^{k+1}-2)\ .$$
Substituting $k=n-\ell-1$, we get
$$(2^{\ell+1}-2) \cdot (2^{\ell+1}-3)\cdot (2^{n-\ell-1}-1) \cdot (2^{n-\ell-1}-2)+(2^\ell-1)^2\cdot (2^{n-\ell}-2)\ .$$

We now sum over the relevant values of $\ell$ and $k$. For case (1), we have
\begin{align*}
&\sum_{\ell=1}^{n-3}\sum_{k=1}^{n-\ell-2} (2^{\ell+1}-2) \cdot (2^{\ell+1}-3) \cdot (2^{k+1}-2)\cdot (2^{k+1}-3) \cdot 2^{3(n-\ell-k-2)})\\
=&\ 5\cdot {2}^{2n}-{\frac {764}{49}}+{\frac {95}{9}}\cdot {2}^{n}-{\frac {25}{3}}
\cdot n\cdot {2}^{n}-n\cdot{2}^{2n}-{\frac {36}{7}}\cdot n+{\frac {16}{441}}\cdot {2}^{3n}\  .
\end{align*}
While in case (2), the number of quartets is
\begin{align*}
&\sum_{\ell=1}^{n-2}\left((2^{\ell+1}-2) \cdot (2^{\ell+1}-3) \cdot (2^{n-\ell-1}-1) \cdot (2^{n-\ell-1}-2) + {(2^{\ell}-1)}^2 \cdot (2^{n-\ell}-2)
\right)\\
=&\ 28+10\cdot n-22\cdot {2}^{n}-6\cdot {2}^{2n}+n\cdot{2}^{2n}+13\cdot n\cdot {2}^{n}\ .
\end{align*}
Summing up the expressions for (1) and (2), the number of ordered quartets with $x_0=0^n$
in $P_{0,1}\cap S_{0,1}$ is
$$
\ {\frac {16}{441}}\cdot {2}^{3n}-2^{2n}+\frac{14}{3}\cdot n\cdot 2^n-\frac {103}{9}\cdot{2}^{n}+
\frac{34}{7}\cdot n +\frac{608}{49} \ .$$
Note that the $\theta(n\cdot2^{2n})$ terms were again cancelled.

\subsection{$P_{0,1}\cap P_{2,3} \cap S_{0,1}$}
We denote the length of the longest common prefix of $x_0,x_1$ by $\ell\  (\ell\leq n-1)$, and the length of the longest common suffix of  $x_0,x_1$ by $k\  (k\leq n-1)$. For $P_{0,1}$, $\ell\geq 1$ should hold, and for $S_{0,1}$, $k\geq 1$ should hold.
Note that in this case, the locations of the longest common suffix and the longest common prefix cannot intersect .
We treat separately the following two cases:
\begin{enumerate}
\item[(1)] $\ell+k+2\leq n$  (the $\ell$ long prefix plus one bit buffer zone, and the $k$ long suffixes plus one bit buffer
zone, do not overlap). Since $1\leq k$, $\ell$ is bounded by $\ell\leq n-3$.
\item[(2)] $\ell+k+1 = n$.
\end{enumerate}
Note that $\ell + k < n$ must hold, for otherwise we would have $x_0=x_1$.
For case (1), by following an argument very similar to case (1) of $P_{0,1} \cap S_{0,1}$, we get that the number
of ordered quartets is
$$2^{3(n-k-\ell-2)}\cdot (2^{\ell+1}-2) \cdot (2^{k+1}-2)\cdot (2^{k+1}-3)\ .$$
Summing over all values of $\ell$ and $k$, we get
\begin{align*}
&\sum_{\ell=1}^{n-3}\sum_{k=1}^{n-l-2} (2^{\ell+1}-2) \cdot (2^{k+1}-2)\cdot (2^{k+1}-3) \cdot 2^{3(n-\ell-k-2)})\\
=&\ {\frac {4}{441}}\,{2}^{3n}-
\frac{1}{3}\,{2}^{2n}+\frac{5}{3}\,n{2}^{n}-{\frac {37}{9}}\,{2}^{n}+{\frac {12}{7}}\,n +{\frac {652}{147}}
\end{align*}

For case (2), given that $x_0$ is $0^n$, we have $x_1=0^\ell 10^k$. The $\ell+1$ prefix of $x_2$ must differ from $0^{\ell+1}$ and
from $0^{\ell}1$, thus there are $2^{\ell+1}-2$ possibilities. Since the longest common prefix of $x_2,x_3$ is also of length $\ell$, the
$\ell+1$ long prefix of $x_2$  determines the $\ell+1$ long prefix of $x_3$. The $k$ long suffix of $x_2$ and of $x_3$ must differ from $0^k$. There are no further constraints, and in particular these two suffixes can be the same, as the next bit of $x_2$ already  differs from that of $x_3$. The number of possibilities for the $k$ long suffix of $x_2$ and of $x_3$ is thus $(2^k-1)^2$. Substituting
$k=n-\ell-1$, the number of ordered quartets in case (2) is
$$(2^{\ell+1}-2)\cdot (2^k-1)^2=(2^{\ell+1}-2)\cdot (2^{n-\ell-1}-1)^2\ .$$
Summing over all values of $\ell$, we get
\begin{align*}
&\sum_{\ell=1}^{n-2}(2^{\ell+1}-2)\cdot (2^{n-\ell-1}-1)^2\\
=&\ \frac{1}{3}\,{2}^{2n}-2\cdot n {2}^{n}+5\cdot{2}^{n}-2\cdot n-\frac{16}{3}\ .
\end{align*}
Adding the two expressions together, we conclude that number of ordered quartets with $x_0=0^n$ in $P_{0,1}\cap P_{2,3}\cap S_{0,1}$ equals
$${\frac {4}{441}}\cdot{2}^{3n}-\frac{1}{3}\cdot n 2^n +\frac{8}{9}\cdot 2^n-\frac{2}{7}\cdot n-\frac{44}{49}\ .$$

\subsection{$P_{0,1}\cap P_{2,3} \cap S_{0,1}\cap S_{2,3}$}
We denote the length of the longest common prefix of $x_0,x_1$ by $\ell\  (\ell\leq n-1)$, and the length of the longest common suffix of  $x_0,x_1$ by $k\  (k\leq n-1)$. Like before, $k,\ell\geq 1$, and we treat separately the following two cases:
\begin{enumerate}
\item[(1)] $\ell+k+2\leq n$  (the $\ell$ long prefix plus one bit buffer zone, and the $k$ long suffix plus one bit buffer
zone, do not overlap). Since $1\leq k$, $\ell$ is bounded by $\ell\leq n-3$.
\item[(2)] $\ell+k+1 = n$.
\end{enumerate}
In case (1), it is (now) easy to see that the number of ordered quartets is
$$2^{3(n-k-\ell-2)}\cdot (2^{\ell+1}-2) \cdot (2^{k+1}-2)\ .$$
Summing over all values of $\ell$ and $k$, we get
\begin{align*}
&\sum_{\ell=1}^{n-3}\sum_{k=1}^{n-l-2} (2^{\ell+1}-2) \cdot (2^{k+1}-2) \cdot 2^{3(n-\ell-k-2)})\\
=&\frac {1}{441}\,{2}^{3n}-\frac{1}{3}\,n{2}^{n}+{\frac {11}{9}}\,{2}^{n}-\frac{4}{7}\,n-{\frac {60}{49}}
\end{align*}
While in case (2), the number of possibilities is
$$(2^{\ell+1}-2)\cdot(2^k-1)=(2^{\ell+1}-2)\cdot(2^{n-\ell-1}-1)\ .$$
Summing over all values of $\ell$, we get
$$\sum_{\ell=1}^{n-2}(2^{\ell+1}-2)\cdot (2^{n-\ell-1}-1)=\ n{2}^{n}-4\cdot{2}^{n}+2\cdot n+4\ .$$
Summing the expressions for case (1) and case (2), we conclude that overall, the number of ordered quartets with $x_0=0^n$ satisfying
$P_{0,1}\cap P_{2,3} \cap S_{0,1}\cap S_{2,3}$ is
$$\frac {1}{441}\,{2}^{3n}+\frac{2}{3}\cdot n 2^n -\frac{25}{9}\cdot 2^n+\frac{10}{7}\cdot n+\frac{136}{49} \ .$$

\section{Putting Everything Together}
Consider the event
$$A= \left(P_{0,1}\cup P_{2,3}\right) \cap \left(S_{0,1}\cup S_{2,3}\right)\ ,$$
A simple manipulation yields
\begin{eqnarray*}
A=    & \left(P_{0,1}\cup P_{2,3}\right) \cap \left(S_{0,1}\cup S_{2,3}\right)\\
    =&  \left(P_{0,1} \cap S_{0,1}\right) \cup \left(P_{0,1} \cap S_{2,3}\right) \cup \left(P_{2,3} \cap S_{0,1}\right) \cup
    \left(P_{2,3} \cap S_{2,3}\right)
\end{eqnarray*}
By the inclusion exclusion principle
\begin{equation*}
\begin{aligned}
\left|A\right|=& \left|\left(P_{0,1} \cap S_{0,1}\right) \cup \left(P_{0,1} \cap S_{2,3}\right) \cup \left(P_{2,3} \cap S_{0,1}\right) \cup\left(P_{2,3} \cap S_{2,3}\right)\right|   \\
=& \left| P_{0,1} \cap S_{0,1}\right|  +
\left|P_{0,1} \cap S_{2,3}\right|  + \left|P_{2,3} \cap S_{0,1}\right| + \left|P_{2,3} \cap S_{2,3}\right|  \\
&-\left|P_{0,1} \cap S_{0,1} \cap S_{2,3}\right| -\left| P_{0,1} \cap S_{0,1}\cap P_{2,3}\right| \\
&-2\left| P_{0,1} \cap S_{0,1}\cap P_{2,3}\cap S_{2,3}\right| -\left|P_{0,1} \cap S_{2,3}\cap P_{2,3}\right| \\
&-\left|P_{2,3} \cap S_{0,1}\cap S_{2,3}\right|
+4\left| P_{0,1} \cap S_{0,1}\cap P_{2,3}\cap S_{2,3}\right| \\
&-\left| P_{0,1} \cap S_{0,1}\cap P_{2,3}\cap S_{2,3}\right| \\
=&\left| P_{0,1} \cap S_{0,1}\right|  +
\left|P_{0,1} \cap S_{2,3}\right|  + \left|P_{2,3} \cap S_{0,1}\right| + \left|P_{2,3} \cap S_{2,3}\right|  \\
&- \left|P_{0,1} \cap S_{0,1} \cap S_{2,3}\right| -\left|P_{0,1} \cap S_{0,1}\cap P_{2,3}\right| \\
&-\left|P_{0,1} \cap S_{2,3}\cap P_{2,3}\right|
-\left|P_{2,3} \cap S_{0,1}\cap S_{2,3}\right| \\
&+\left| P_{0,1} \cap S_{0,1}\cap P_{2,3}\cap S_{2,3}\right| \\
=&2\left| P_{0,1} \cap S_{0,1}\right|  +
2\left|P_{0,1} \cap S_{2,3}\right|-4\left|P_{2,3} \cap S_{0,1}\cap S_{2,3}\right|+\left| P_{0,1} \cap S_{0,1}\cap P_{2,3}\cap S_{2,3}\right|
\end{aligned}
\end{equation*}
Substituting the expressions we derived for the various subsets, we conclude that the
number of ordered quartets with $x_0=0^n$ in $A$ equals
$$\frac{1}{9}\,{2}^{3n}-2\cdot2^{2n}+{\frac {20}{3}}\,n{2}^{n}-\frac{127}{9}\,{2}^{n}+6\,n+16\ .$$
Removing the $x_0=0^n$ restriction, the number of ordered quartets in $A$ equals
$$\frac{1}{9}\,{2}^{4n}-2\cdot 2^{3n}+{\frac {20}{3}}\,n{2}^{2n}-\frac{127}{9}\,{2}^{2n}+6\,n{2}^{n}+16\cdot{2}^{n}\ .$$
We now introduce two related sets, $B$ and $C$:
$$B= \left(P_{0,2}\cup P_{1,3}\right) \cap \left(S_{0,2}\cup S_{1,3}\right)\ , C= \left(P_{0,3}\cup P_{1,2}\right) \cap \left(S_{0,3}\cup S_{1,2}\right)\ .$$
Clearly $A, B, C$ are mutually exclusive and $A, B, C$ have the same number of ordered quartets.
Therefore $$\left| A \cup B \cup C\right| = \frac{1}{3}\,{2}^{4n}-6\cdot 2^{3n}+20\cdot n{2}^{2n}-\frac{127}{3}\,{2}^{2n}+18\cdot n{2}^{n}+48\cdot{2}^{n}\ .$$
We observe that the union $A \cup B \cup C$ contains exactly those ordered quartets
on $x_0, x_1, x_2, x_3$ that agree in both prefix and suffix trees.

\section{Unordered Quartet and the Quartet Distance}
So far, we counted ordered quartets.
In the quartet distance problem, we are interested in {\it unordered} quartets and not in ordered ones.
There are $4!=24$ permutations over a set of 4 distinct elements, $\{x_0,x_1,x_2,x_3\}$. We will show that
for any set $\{x_0,x_1,x_2,x_3\}$, either the suffix and the prefix tree agree
for all the ordered 4-tuples corresponding to these 24 permutations,
(namely the ordered event $A \cup B \cup C$, is satisfied), or there is
no agreement for any of the permutations. This statement implies that the number of unordered quartets
where the two trees agree is exactly this number for the ordered case, divided by 24.

We will show that $A$ is invariant under exactly 8 permutations of ordered $4$-tuples.
A different set of
8 permutations maps ordered $4$-tuples that satisfy $A$ to different orders where the
$4$-tuple satisfies $B$, and yet another
8 permutations map ordered $4$-tuples that satisfy $A$ to different orders satisfying $C$.

Suppose the ordered pair $(x_0,x_1,x_2,x_3)$ satisfies $A$, namely
$$\left(P_{0,1}(x_0,x_1,x_2,x_3)\cup P_{2,3}(x_0,x_1,x_2,x_3)\right) \ \cap \  \left(S_{0,1}(x_0,x_1,x_2,x_3)\cup S_{2,3}(x_0,x_1,x_2,x_3)\right)\ .$$
Membership in $A$ is invariant under each of the following 3 permutations
and their compositions:  Transposing $x_0, x_1$; transposing $x_2, x_3$;
replacing $x_0, x_1$ by $x_2, x_3$. These 3 permutations generate a subgroup of
size $8$.

Starting with an ordered quartet $(x_0,x_1,x_2,x_3)$ in $A$, and transposing $x_1$ with $x_2$,
the new ordered quartet  $(x_0,x_2,x_1,x_3)$  is now in $B$. By first applying one
 of the 8 permutations keeping $(x_0,x_1,x_2,x_3)$ in $A$, and then this
transposition, we conclude that there is a coset of 8 permutations, moving an ordered quartet from $A$
to $B$. A similar argument holds regarding moving from $A$ to $C$, employing the transposition
of $x_1$ with $x_3$.  Clearly, the same argument is applicable if we  start with an ordered
quartet that satisfies $B$ or $C$.

We conclude that if the prefix and suffix trees agree on one ordered quartet, then
they will agree
on all 24 permutations of it. Dividing the number of ordered permutations on $n$-bit strings where
this event occurs by $24$, we conclude that the number of {\em unordered} permutations
that agree equals
$$ \left(\frac{1}{3}\,{2}^{4n}-6\cdot 2^{3n}+20\cdot n{2}^{2n}-\frac{127}{3}\,{2}^{2n}+18\cdot n{2}^{n}+48\cdot{2}^{n}\right)/24\ .$$
The number of unordered quartets equals
$$\frac{2^n\cdot(2^n-1)\cdot(2^n-2)\cdot(2^n-3)}{24}=
\frac{2^{4n}-6\cdot 2^{3n}+11\cdot 2^{2n}-6\cdot 2^n}{24}\ .$$

Thus the quartet distance between the two trees equals
$$\frac{2^{4n}-6\cdot 2^{3n}+11\cdot 2^{2n}-6\cdot 2^n - \left(\frac{1}{3}\,{2}^{4n}-6\cdot 2^{3n}+20\cdot n{2}^{2n}-\frac{127}{3}\,{2}^{2n}+18\cdot n{2}^{n}+48\cdot{2}^{n}\right)}{24} $$
$$= \frac{\frac{2}{3}\cdot 2^{4n}-20\cdot n2^{2n}+\frac{160}{3}\cdot 2^{2n}-18\cdot n2^n-54\cdot2^n}{24} \ .$$

The ratio, or {\em normalized} quartet distance for the suffix and prefix trees on $N=2^n$ leaves
equals
$$\frac{\frac{2}{3}\cdot 2^{4n}-20\cdot n2^{2n}+\frac{160}{3}\cdot 2^{2n}-18\cdot n2^n-54\cdot2^n}
{2^{4n}-6\cdot 2^{3n}+11\cdot 2^{2n}-6\cdot 2^n } \ .$$

It is easy to see that this ratio indeed converges to $2/3$ as $n\to\infty$. What is not so obvious
is that this ratio is a monotonically decreasing function of $n$.
For small values of $n$ we get the following distances and ratios. For these values (and many others
we tested numerically) the ratio indeed decreases monotonically with growing values of $n$.
$$
\left|\begin{array}{c|c|c|c|c|c|c|c|c|} \hline
n & 3 & 4 & 5 & 6 & 7 & 8 & 9 & 10  \\
\hline \mbox{distance} & 60 & 1452 & 26944 & 454224 & 7396416 & 119011264 & 1907486208 & 30535571712  \\
\hline \mbox{ratio} & 0.857 & 0.797 & 0.749 & 0.714 & 0.693 & 0.680 & 0.674 & 0.670\\ \hline \end{array}\right.
$$
To prove the above mentioned monotonicity, we show that the ratio for $N=2^n$ is larger than the
ratio for $N_1=2^{n+1}$. The proof involves somewhat tedious yet elementary
arithmetic manipulations.

Let $R(n)$ denote the ratio (number of disagreeing unordered quartets, divided by the
total number of unordered quartets) for the prefix and suffix trees with $N=2^n$ leaves. Then
$$    R(n)= \frac{\frac{2}{3}\cdot 2^{4n}-20\cdot n2^{2n}+\frac{160}{3}\cdot 2^{2n}-18\cdot
    n2^n-54\cdot2^n}{2^{4n}-6\cdot 2^{3n}+11\cdot 2^{2n}-6\cdot 2^n }  $$
$$ R(n+1)=   \frac{\frac{2}{3}\cdot 2^{4(n+1)}-20\cdot (n+1)2^{2(n+1)}+\frac{160}{3}\cdot 2^{2(n+1)}-18\cdot (n+1)2^{n+1}-54\cdot2^{n+1}}
{2^{4(n+1)}-6\cdot 2^{3(n+1)}+11\cdot 2^{2(n+1)}-6\cdot 2^{n+1} } \ .$$
To show that $R(n)>R(n+1)$, we first compute
$$\mbox{numerator}(R(n))\cdot \mbox{denominator}(R(n+1)) -
\mbox{numerator}(R(n+1))\cdot \mbox{denominator}(R(n))\ ,$$ and then take the derivative of this expression with respect to the (real) variable $n$. The result (obtained using Maple) equals
\begin{align*}
-&3792\cdot{2}^{4n}\ln  \left( 2 \right) +1440\cdot2^{4n}\cdot n \cdot \ln  \left( 2 \right)+360\cdot2^{4n} -240\cdot\ln  \left( 2 \right) 2^{5n}\\
-& 342\cdot2^{3n}-1026 \cdot2^{3n}\cdot n \cdot \ln  \left( 2 \right) +10908\cdot2^{3n}\ln  \left( 2 \right)\\
-&1944\cdot2^{2n}\cdot n \cdot \ln  \left( 2 \right)-972\cdot2^{2n}-7824\cdot2^{2n}\ln  \left( 2 \right)\\
+& 954\cdot{2}^{n}\cdot n \cdot \ln  \left( 2 \right)+954\cdot{2}^{n}
+ 948\cdot{2}^{n}\ln  \left( 2 \right)
\end{align*}
It is easy to see that for $n\geq 11$, the derivative is negative, as following:
The first term in the first line, $-3792\cdot{2}^{4n}\ln  \left( 2 \right)$, dominates the third term in the same line. The difference $1440\cdot2^{4n}\cdot n \cdot \ln  \left( 2 \right)-240\cdot\ln  \left( 2 \right) 2^{5n}$ is negative for all $n\geq 5$.
For $n\geq 11$, the second term in the second line, $-1026 \cdot2^{3n}\cdot n \cdot \ln  \left( 2 \right)$,
dominates the third term in the same line, $+10908\cdot2^{3n}\ln  \left( 2 \right)$.
Each of the three terms containing $2^n$ (fourth line) is dominated by a term containing
$2^{2n}$ (third line) with a minus sign. Finally,
for (integer) values of $n$ in the range $3\leq n \leq 11$, direct computation verifies that $R(n)-R(n+1)>0$.

\section{Concluding Remarks and Open Problems}
\label{conc}
There is more than a single notion of what an ``explicit construction'' means. Possibly the most popular one is that an explicit construction is (1) deterministic, and (2) it runs in polynomial time (polynomial in the size of the object being constructed). Under this definition, by de-randomizing a randomized labeling of the leaves, we would get an explicit construction with quartet distance being asymptotically $\frac{2}{3}\binom{N}{4}$. It may require some additional work to determine by how much the exact bound resulting from this approach exceeds $\frac{2}{3}\binom{N}{4}$ for concrete values of $N$. This construction is deterministic, and its running time is polynomial in the size of the resulting trees, $N=2^n$. Thus, this is an explicit construction by the definition above. Furthermore, it is applicable to any two trees (not just complete, balanced binary trees), and any size $N$ (not just a power of 2). On the other hand, it is hard to argue that (for complete, balanced binary trees) our construction is much simpler, and arguably elegant, than what the de-randomization yields.

A ``strongly explicit construction'' enables one to determine, given the specification of an entry in the object, the
contents of this entry, in time that is polynomial in the length of the description of the entry (as opposed to the size
of the complete object). This is applicable to a variety of objects, {\em e.g.} graphs, matrices, and codes \cite{Indyk}.
In our context, a strongly explicit construction should be able to determine, in time polynomial in $n$ (and {\em not} in $2^n$)
the label of a leaf, given the description of this leaf. Furthermore, given the labels of four leaves, we should be able to determine
the induced quartets topologies for the two trees.

The standard de-randomized construction is {\em not} strongly explicit. Essentially, it mimics the randomized
construction, where one first  assigns labels to all leaves, and only then can
determine the labels of specific leaves or the topologies of specific quartets. By way of contrast, our prefix--suffix construction
is strongly explicit. Assuming we use the standard labeling of the complete, balanced binary trees by the prefix order, then the labeling of, say,
{\tt 0111} in the prefix tree will be, well, {\tt 0111}, which is the rightmost leaf on the major left subtree.
In the suffix tree it will be placed on the one left to the rightmost leaf, the one labeled by {\tt 1110} in prefix order (we simply reverse
the binary string to move from prefix to suffix order). So determining the location is done in linear time, using a trivially simple algorithm.

Turning to quartets given four labels, in the prefix order the two labels with longest common prefix will be together, and dually for
the suffix order. So determining prefix and suffix quartets topologies is also done, given the four labels, by a trivial linear time algorithm.
See the following figure, for labels of length $n=4$.

\bigskip
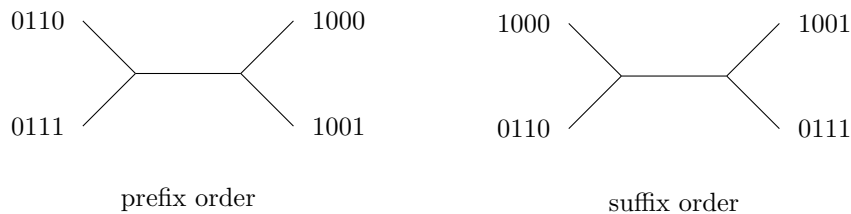
\begin{figure}[h]
\centering
\begin{tikzpicture}[scale=1.4]
\draw  (0,0) -- (1,0);
\draw  (0,0) -- (-.5,-.5);
\draw  (0,0) -- (-.5,.5);
\draw  (1,0) -- (1.5,.5);
\draw  (1,0) -- (1.5,-.5);
\node  at (-.5,-.5) [label=left:$0111$] {};
\node  at (-.5,.5) [label=left:$0110$] {};
\node  at (1.5,-.5) [label=right:$1001$] {};
\node  at (1.5,.5) [label=right:$1000$] {};
\node  at (.5,-1.2) {prefix order};
\end{tikzpicture}
\qquad\qquad
\begin{tikzpicture}[scale=1.4]
\draw  (0,0) -- (1,0);
\draw  (0,0) -- (-.5,-.5);
\draw  (0,0) -- (-.5,.5);
\draw  (1,0) -- (1.5,.5);
\draw  (1,0) -- (1.5,-.5);
\node  at (-.5,-.5) [label=left:$0110$] {};
\node  at (-.5,.5) [label=left:$1000$] {};
\node  at (1.5,-.5) [label=right:$0111$] {};
\node  at (1.5,.5) [label=right:$1001$] {};
\node  at (.5,-1.2) {suffix order};
\end{tikzpicture}
\end{figure}

As noted above, our construction and proof are applicable only to complete, balanced binary trees on $N=2^n$  leaves.
It will be interesting to extend these results to other tree topologies, and also values of $N$ that are not exact power of $2$. 
We note that the tree topology may have a substantial impact on the feasibility of a proof. For example, 
Alon, Naves, and Sudakov \cite{ANS} have shown a $\left(0.69 +  o(1)\right) \cdot\binom{N}{4}$ upper bound on the normalized quartet distance of general binary trees (for large enough $N$), but a better $\left(2/3 +  o(1)\right) \cdot\binom{N}{4}$ upper bound for caterpillar trees. 
Finally, it will 
be interesting to prove or refute the conjecture that for large enough $n$, the largest quartet distance
on trees with $N=2^n$ leaves is obtained by the suffix and prefix trees.

\section*{Acknowledgements}Thanks to Stefan Gr\"unewald for helpful discussions, which motivated this work. We would also
like to thank Noga Alon and Gil Cohen for their help regarding what ``explicit construction'' exactly means. PLE was supported in
part by Hungarian NSF grant K116769. Part of this work was done when PLE visited BC, supported by an exchange program of the Hungarian and Israeli academies of Science. BC was supported by a grant from the Blavatnik Computer Science Research Fund, and by the LTZI (Long Term Zero Income) fund of the ISF (Israeli Science Foundation).

\bibliographystyle{abbrv}

\end{document}